\documentclass[letterpaper, 10 pt, conference]{ieeeconf}  

\IEEEoverridecommandlockouts                              
\usepackage{dsfont,soul,color,url}
\usepackage{float}
\usepackage{multicol}
\usepackage{graphicx}
\usepackage{subcaption}
\usepackage[algo2e,ruled,lined]{algorithm2e}
\usepackage{amsmath, amssymb, amsfonts}
{
      \newtheorem{corollary}{Corollary}
      \newtheorem{theorem}{Theorem}
      
      \newtheorem{lemma}{Lemma}
      
}

\newcommand{\norm}[1]{\left\lVert#1\right\rVert}

\title{\LARGE \bf
Input design for the optimal control of networked moments
}

\author{Philip Solimine and Anke Meyer-Baese
\thanks{This work was not supported by any organization}
\thanks{P. Solimine is with the Vancouver School of Economics at the University of British Columbia, Vancouver, BC, Canada
        {\tt\small philip.solimine@ubc.ca}}%
\thanks{A. Meyer-Baese is with the Department of Scientific Computing at the Florida State University, Tallahassee, FL, USA
        {\tt\small ameyerbaese@fsu.edu}}%
}

\begin{document}

\maketitle
\thispagestyle{empty}
\pagestyle{empty}

\begin{abstract}
We study the optimal control of the mean and variance of the network state vector. We develop an algorithm that uses projected gradient descent to optimize the control input placement, subject to constraints on the state that must be achieved at a given time threshold; seeking to design an input that moves the moment at minimum cost. First, we solve the state-selection problem for a number of variants of the first and second moment, and find solutions related to the eigenvalues of the systems' Gramian matrices. We then nest this state selection into projected gradient descent to design optimal inputs.
\end{abstract}

\captionsetup{width=\linewidth}

\section{Introduction}
Finding a set of leaders, drivers, or actuators from which the dynamics of the system can be controlled most efficiently is an important problem, and one which has recently received an influx of attention \cite{li2015minimum,gao2020optimal,lindmark2018minimum,jia2013control,pasqualetti2014controllability,lindmark2020centrality,baggio2019gramian,klickstein2020selecting}. We approach the problem from a new angle; rather than taking the desired final state as given, we search first for the {\it state} that satisfies the goal at minimum cost. Then, we examine the placement of inputs to optimally adjust the output measure.

An example of a network system with linear dynamics is the classic DeGroot model of opinion dynamics \cite{degroot1974reaching}. In this system, the state of the system represents the opinions of a set of nodes, connected by a social network whose structure determines the dynamics matrix $A$. The control of such systems is thus a topic of great interest to social scientists \cite{mostagir2019society,liu2014control,acemoglu2010spread}. Thus it is natural to consider how input nodes could be placed efficiently in order to, for example, change the average opinion over some topic \cite{vohra2020strategic}, induce discord in opinions \cite{gaitonde2020adversarial,galeotti2021discord}, or encourage fast convergence \cite{golub2012homophily}. In particular, recent work has used linear control systems and continuous-time DeGroot-type settings to \cite{candogan2022social} to model the influence of digital platforms on social behavior. In this way, we provide a set of results to relate the emerging literature on platform design for social influence with the rich and growing space of results for the optimal control of network systems.

Previous literature on input placement for efficient target control of network systems \cite{gao2014target,li2019target,klickstein2020selecting}, has focused on fully controlling a subset of nodes, rather than sufficient statistics of the state. Much literature on network control input placement has also focused on driver node selection, as opposed to augmentation, in which each input influences a set of existing nodes to different degrees \cite{lindmark2020centrality,gao2020optimal}. A notable exception is \cite{li2015minimum}, which develops the Projected Gradient Method (PGM) by restricting the columns of the control input matrix, representing its schematic of outgoing connections, to be embedded on a sphere surface rather than restricted to be drawn from the canonical basis.

We frame these problems of input design as optimal control for networks with state distribution goals. These goals come in the form of statistical features of a distribution. Specifically, the mean and the second moment or variance. With mean constraints, the problem is straightforward and can be solved using standard methods. The variance or second moment problem can be solved by nesting convex programming in the general case, and allows for a closed-form bound on the energy cost derived from the spectrum of the controllability Gramian.

\section{Preliminaries}
\subsection{Framework}
We focus on the case of linear dynamical systems, interpreted as weighted digraphs. The system will, in general, be referred to by its dynamics matrix $A$. We will define the graph $\mathcal{G}$ as a tuple consisting of a set of $n$ vertices (nodes) $V$ (indexed by the positive integers $\mathbb{Z}_+$), and a set of directed edges (links) $E \subset V\times V$ as ordered pairs of nodes. We say that an edge $(i,j) \in E$ if and only if $A_{ij} \neq 0$. Each edge $(i,j)$ may also be associated with a weight equal to its corresponding entry in the dynamics matrix $w_{ij} = A_{ij}$. In addition, each node in $V$ possesses a state variable $x_i \in \mathbb{R}$, collected into the state vector $x \in \mathbb{R}^n$. Further, we focus on the case of a linear, time-invariant system in which $A$ is fixed over time.

\subsection{Dynamics and control of linear systems}
We focus on the case of finite-time control of a canonical linear time-invariant system with the following dynamics:
\begin{equation}\label{dynamics}
    \dot{x} = Ax + Bu,
\end{equation}
where $x$ is a $n\times 1$ column vector containing the states of $n$ variables, $\dot{x}\triangleq\frac{\partial x}{\partial t}$, the $n\times n$ matrix $A$ describes the behavior of the autonomous system, $u$ is a $m\times 1$ vector containing control input values and $B$ is a $n\times m$ matrix which describes the nodes to which input from $u$ is sent. If $B$ is binary, with exactly one $1$ per column, then the nodes corresponding to these $1$'s are referred to as driver nodes. Each column of the control input schematic $B$ is associated with unique input signal, and thus corresponds with an augmented node.

The final state of this system at finite terminal time $t^*$, with no control input, is given by $e^{t^* A}x_0$. For convenience, we will use $z\triangleq e^{t^* A}x_0$ to denote the final state of the autonomous dynamics. Existing results show that this system can be controlled to any state if the Kalman rank condition holds -- if the controllability matrix $C$ is full rank \cite{liu2011controllability}. The controllability matrix is given by
\begin{equation}
    C = [\; B,\; AB,\; A^2B,\; \dots, \; A^{n-1}B\;], \label{kalman}
\end{equation}
where the comma denotes horizontal concatenation. Thus to check the Kalman condition is to check that the following holds:\begin{equation}
    \text{rank}(C) = n.
\end{equation}

Later on, when we turn our focus to designing input schematics, we may switch this notation to $C(B)$ or $C(B;A)$ interchangeably, to highlight that this matrix is functionally dependent on the input placement.
\section{Controlling the average}
\subsection{State selection with average-state goals}
In this section, we will use the standard method of Lagrange multipliers to derive a minimum-energy final state for a system subject to a linear constraint on its central tendency, such as the average value of nodes at evaluation time $t^*$. In this case, constraints take the form\begin{equation}
\gamma^\top x = \eta\label{meanconstraint},
\end{equation}
for some threshold $\eta \in \mathbb{R}$ and vector $\gamma\in \mathbb{R}^n$.

Under the optimal input signals, the minimum energy required to control the system from initial state $x_0$ to final state $x \in \text{span} (C)$ by the time $t=t^*$ (beginning, without loss of generality, at time $t=0$) is known \cite{lindmark2018minimum} to have the following closed-form solution\begin{equation}
    \mathcal{E}_{t^*}(x) = (z - x)^\top W^{\dagger} (z - x)\label{energy},
\end{equation}
where $W^{\dagger}$ is the Moore-Penrose pseudoinverse of the following matrix called the {\it reachability Gramian}\begin{equation}\label{gram}
    W(B;A,t^*) = \int_0^{t^*} e^{tA}BB^\top e^{tA^\top}dt.
\end{equation} The notations $W$, $W(B)$, and $W(B;A,t^*)$ will be used interchangably depending on context. Initially, we will focus on the situation in which $B$ is fixed, in which case it will be convenient to drop the functional notation $W(B)$ in favor of simply $W$. Whenever a $W$ appears, however, it will be important to remember that it may be viewed as a function of the control schematic $B$ in the case of endogenous input placement. The invertibility of $W$ is guaranteed provided that the system is controllable under the given control schematic $B$.

In this section, we will take the control input schematic $B$ to be given, but we will relax this assumption later in order to explain how switching driver sets can lead to savings on control input energy cost. Using this, it is straightforward to derive the final state $x^*$ which satisfies the constraint at the lowest possible cost. That is, the optimal final state is given by the solution of the minimization of (\ref{energy}) subject to constraint (\ref{meanconstraint}). We will use $\mathds{1}$ to denote an $n\times 1$ column vector of ones and define $\alpha(\gamma) \triangleq \gamma^\top e^{t^* A}x_0 - \eta$ for convenience. We will denote this as a constant $\alpha$ because it is independent of $B$. 

For convenience, we will define the scalar-valued function $\kappa :\mathbb{R}^{n\times m}\times\mathbb{R}^n \rightarrow\mathbb{R}$ as\begin{equation}
    \label{kappa}\kappa(B;\gamma) \triangleq \gamma^\top W(B) \gamma.
\end{equation}

Then the solution to the problem of state selection subject to mean constraints gives a simple closed-form solution which will comprise our first theorem:\begin{lemma}\label{thm:mean_state}
 The optimal state which can satisfy (\ref{meanconstraint}) at the lowest possible energy cost for a given control input matrix $B$ and linear output $\gamma$ is given by \begin{equation}x_a^* = z - \frac{\alpha(\gamma)}{\kappa(B,\gamma)} W \gamma. \label{xastar}\end{equation}
\end{lemma}
\begin{proof}This is the solution of the optimization problem to minimize (\ref{energy}), subject to the constraint (\ref{meanconstraint}), and the constraint that $x=W^\dagger W x$ (optimizing over the space of reachable states). The Lagrangian for this optimization problem is \begin{equation}
    \mathcal{L} = (e^{t^* A}x_0 - x)^\top W^{\dagger}(e^{t^* A}x_0 - x) + \psi (\gamma^\top x - \eta),
\end{equation}
where the Lagrange multiplier $\psi$ can be interpreted as the shadow price of satisfying the constraint \ref{meanconstraint}. Using the symmetry of W and taking the gradient of $\mathcal{L}$ with respect to the elements of $x$ gives the following set of first-order conditions: $\nabla_{x} \mathcal{L}(x_a^*) = 2W^{\dagger}(x_a^* - e^{t^* A}x_0) + \psi \gamma \triangleq 0$. Rearranging yields \begin{equation} \frac{\psi}{2} W \gamma = e^{t^* A}x_0 - x_a^*.\label{foc:x}\end{equation} The derivative with respect to the Lagrange multiplier $\psi$ is $\frac{\partial \mathcal{L}}{\partial \mathcal{\psi}} = \gamma^\top x^*_a - \eta \triangleq 0$ which recovers the constraint as a first-order condition \begin{equation}\gamma^\top  x^*_a = \eta.\label{foc:rho}\end{equation} Premultiplying (\ref{foc:x}) by $\gamma^\top$ gives \begin{equation}\frac{\psi^*}{2}\gamma^\top W\gamma = \gamma^\top e^{t^* A}x_0 - \gamma^\top x^*_a. \label{prefinal}\end{equation} Substituting in (\ref{foc:rho}) and rearranging gives the closed-form solution for the shadow price:\begin{equation}
    \psi^* = 2\left( \frac{\gamma^\top e^{t^* A}x_0 - \eta}{\gamma^\top W\gamma} \right) = 2 \frac{\alpha}{\gamma^\top W(B) \gamma}.\label{shadowprice}
\end{equation} Finally, substituting (\ref{shadowprice}) back in to (\ref{foc:x}) and solving for $x^*$ gives the above solution for optimal state.\end{proof}

The following corollary immediately follows which gives a closed-form solution to the minimum energy with which the constraint can be satisfied, given a fixed control schematic $B$:
\begin{lemma}\label{cor:mean_energy}
The minimum energy cost required to achieve the single linear constraint goal (\ref{meanconstraint}) with output $\gamma$ under the control schematic $B$ is given by:\begin{equation}
    \bar{\mathcal{E}}_{t^*}^\eta (B,\gamma; A, \eta, x_0) = \frac{\alpha(\gamma;x_0,\eta)^2}{\kappa(B,\gamma)}\label{minenergyavg}.
\end{equation}
\end{lemma}
\begin{proof}This derivation follows by substituting the optimal state $x_a^*$ from (\ref{xastar}) into the minimum energy formula (\ref{energy}). The $e^{t^* A}x_0$ immediately cancel. This leaves:\begin{align*}
    \bar{\mathcal{E}}_{t^*}^\eta (B; A, \eta, x_0) &= \left(\frac{\alpha}{\gamma^\top W \gamma} W \gamma\right)^\top W^{\dagger} \left(\frac{\alpha}{\gamma^\top W \gamma} W \gamma\right)\\
     &= \left(\frac{\alpha}{\gamma^\top W \gamma}\right)^2 (W\gamma)^\top W^{\dagger} (W\gamma)\\
     &= \left(\frac{\alpha}{\gamma^\top W \gamma}\right)^2 (\gamma^\top W^\top) W^{\dagger} (W\gamma).
\end{align*}
$W$ is symmetric so that $W^\top = W$. Finally we have \begin{equation}
    \bar{\mathcal{E}}_{t^*}^\eta (B; A, \eta, x_0) = \frac{\alpha^2}{\gamma^\top W(B) \gamma}.
\end{equation}
Equation (\ref{minenergyavg}) follows directly, by substituting the definition of $\kappa$. \end{proof}

\subsection{Input placement with average-state goals}
In equation (\ref{minenergyavg}) we have derived a closed form expression for the energy cost associated with meeting the mean constraint, conditional on a given input schematic $B$. To optimize the input placement, we will use matrix calculus to derive an analytic solution to the Optimal Minimum Augmentation Placement (OMAP) problem for control of a linear combination of network states, such as the average, and its associated energy cost.

Classical results on {\it output controllability} give us the following lemma.\begin{lemma}
For any linear system with a single linear output $\gamma\in\mathbb{R}^n$, there exists a one-dimensional control input $B\in\mathbb{R}^{n}$ such that the augmented system is output controllable. \label{onecontroller}
\end{lemma}\begin{proof}This follows an established result in output controllability. The lemma be formulated as a special case of the classical {\it output controllability} studied by \cite{morse1971output}. This problem discusses the controllability of an output described as:\begin{equation}
    y = Ox + Zu.
\end{equation}

Then the problem can be viewed as a solution to output controllability for the special case where $Z=0$ and $O$ is one-dimensional\footnote{In this way, our work is closely related to the problem of \textit{target control} \cite{gao2014target,duan2019energy,klickstein2020selecting}, which refers to a different special case of output controllability, in which $D=0$ and the rows of $C$ are drawn from those of an $n\times n$ identity matrix.}. For a full proof, see \cite{morse1971output}.\end{proof}

This special case with a single, mixed output emits a closed-form solution for both the energy cost and the optimal control input schematic. In particular, recent studies on target control become focused on target node selection and subsequent optimization of the control driver or actuator sets \cite{gao2020optimal,klickstein2020selecting,chen2020optimizing}.

Lemma \ref{onecontroller} contains, as a special case, control of the average state: \begin{lemma}
For any network, a single, suitably designed control input is sufficient to fully control the average state to any desired threshold $\eta$.
\end{lemma}\begin{proof}
The average state is simply a linear combination of network states with equal weights, with $\gamma=\frac{1}{n}\mathds{1}$.
\end{proof}

For example, consider the empty network $A_0 = 0_{n\times n}$, with the control schematic $B_i = e_i$ (where $e_i$ denotes a vector of zeroes with a one in the $i_{th}$ position). Clearly this system is not controllable according to the Kalman or PBH conditions. In this network, however, since a single node can be controlled, we can easily change the {\it average} state of the nodes only by changing the state of a single node. Since $\frac{1}{n}x^\top \mathds{1} = \frac{1}{n}\sum_{j=1}^n x_j$, we simply need to move the state of node $i$ to be $x_i = n\eta - \sum_{j=1,j\neq i}x_j$ to compensate for the states of the uncontrolled nodes.

The optimization problem we consider in this section is the following:\begin{equation}\begin{aligned}
    &\underset{B\in \mathbb{R}^{n\times m}}{\text{minimize}} & &\bar{\mathcal{E}}_{t^*}^\eta(B,\gamma)\\&\text{subject to}&&tr(B^\top B) = m
    \end{aligned}\label{meaninputselection}
\end{equation}
That is, our goal is to find a control schematic $B_a^*$ which minimizes the energy required to meet the mean constraint (\ref{meanconstraint}). The constraint in (\ref{meaninputselection}) is simply a normalization condition on the control input matrix, of the type used in \cite{li2015minimum}. The solution to (\ref{meaninputselection}) gives our next theorem.

In particular, we will relate the solutions of the optimization problem to the spectrum of a certain matrix $\Phi(\gamma)$, which we will define as:\begin{equation}
    \Phi = \int_0^{t^*} e^{t A^\top} \gamma\gamma^\top e^{t A} dt.
\end{equation}
\begin{theorem}\label{optimalcontrol}
The columns of the optimal control augmentation $B_a^*$ that solve the optimization problem (\ref{meaninputselection}) are eigenvectors associated with the largest eigenvalue of $\Phi(\gamma)$.
\end{theorem}
\begin{proof}
First, note that for all conforming matrices $X$, $Y$, $Z$, and $Q$, we have:\begin{equation}\label{difflemma}
\frac{\partial}{\partial X} tr(YXQX^\top Z ) = ZYXQ + Y^\top Z^\top X Q^\top.\end{equation}
We can solve the optimization problem using the method of Lagrange multipliers and techniques from matrix calculus. First, we can use that $\bar{\mathcal{E}}_{t^*}^\eta = \frac{\alpha^2}{\kappa (B,\gamma)}$ to reduce the optimization problem to a simpler one. That is, we note that the minimum energy is inversely proportional to the sum of elements in the controllability Gramian $\kappa(B,\gamma)$. since $\alpha$ is constant with respect to $B$, the $B_a^*$ which minimizes the energy is thus given by the one which maximizes the denominator, $\kappa(B,\gamma) = \gamma^\top W(B) \gamma$. Thus the new optimization problem is:\begin{equation}\begin{aligned}
    &\underset{B\in \mathbb{R}^{n\times m}}{\text{maximize}} & &\kappa(B,\gamma)\\&\text{subject to}&&tr(B^\top B) = m
    \end{aligned}\label{meaninputselection2}
\end{equation}
First, note that the trace of a scalar is equal to that scalar. This gives $\kappa(B) = \gamma^\top W(B) \gamma = tr(\gamma^\top W(B) \gamma)$. The Lagrangian for this optimization problem may then be written as:\begin{equation}
    \mathcal{L} = tr(\gamma^\top W(B) \gamma) - \psi (tr(B^\top B) - m).
\end{equation} Taking the first-order condition with respect to the Lagrange dual variable $\psi$ recovers the constraint\begin{equation}
    tr(B^\top B) = m.
\end{equation}
The first order condition with respect to $B$ requires matrix calculus to derive. Note that the condition can be written as:\begin{equation}
    \label{focb} \frac{\partial \mathcal{L}}{\partial{B}} = \frac{\partial}{\partial B}\kappa (B) - \psi \frac{\partial}{\partial B} tr(B^\top B) \triangleq 0.
\end{equation}
Since $tr(B^\top B) = \sum_{i=1}^n \sum_{j=1}^m B_{ij}^2$, it follows that $\frac{\partial}{\partial B} tr(B^\top B) = 2B$. $\frac{\partial}{\partial B} \kappa(B)$ can be derived as follows. Note that we have:
\begin{equation}
\begin{aligned}
    \kappa(B) &= \gamma^\top W(B) \gamma = tr(\gamma^\top W(B) \gamma)\\
    &= tr\left( \gamma^\top \int_0^{t^*} e^{t A}B B^\top e^{t A^\top} dt \gamma \right)\\
    &=tr\left( \int_0^{t^*} \gamma^\top e^{t A}B B^\top e^{t A^\top} \gamma dt \right)\\
    &=\int_0^{t^*} tr(\gamma^\top e^{t A}B B^\top e^{t A^\top} \gamma)dt.
\end{aligned}\label{reduction}\end{equation}

We can then use the Leibniz integral rule to interchange integration with differentiation. Thus we have:\begin{equation}
    \begin{aligned}
    \frac{\partial}{\partial B} \kappa(B) &= \frac{\partial}{\partial B} \int_0^{t^*} tr(\gamma^\top e^{t A}B B^\top e^{t A^\top} \gamma)dt\\
    &= \int_0^{t^*} \frac{\partial}{\partial B} tr(\gamma^\top e^{t A}B B^\top e^{t A^\top} \gamma)dt.
    \end{aligned}
\end{equation}
From here, we can use (\ref{difflemma}), (substituting $X\triangleq B$, $Q \triangleq I$, and $Y \triangleq Z^\top \triangleq \gamma^\top e^{t A}$), which yields the following:\begin{equation}
    \frac{\partial}{\partial B}\kappa(B) = 2\int_0^{t^*} e^{t A^\top} \gamma\gamma^\top e^{t A^\top} dt B = 2\Phi B.
\end{equation}
Plugging this back into the general first-order condition for $B$ (\ref{focb}) and rearranging gives:

\begin{equation}
    \Phi B_a^* = \psi^* B_a^*\label{focb2}.
\end{equation}

It follows that the columns of $B_a^*$ must each be eigenvectors of $\Phi$ associated with the {\it same eigenvalue}. This means that the local optima of the energy function, along the sphere surface $tr(B^\top B) = m$, are associated with an eigenvalue of $\Phi$. 

Denote this eigenvalue as $\lambda^* \triangleq \psi^*$ to highlight that the local optima of the Lagrangian are determined by the eigenvalues of $\Phi$. What remains is to show that $\lambda^*$ is the top (largest) eigenvalue. To do this, recall the formula for minimum energy (\ref{minenergyavg}). Again, this is that $\bar{\mathcal{E}}^\eta_{t^*}(B) = \frac{\alpha^2}{tr(\gamma^\top W(B) \gamma)}$. Our objective was to maximize the denominator of this energy cost, which in (\ref{reduction}) we saw could be written as: $\kappa(B) = \int_0^{t^*} tr(\gamma^\top e^{t A}B B^\top e^{t A^\top} \gamma)dt$. Using properties of the trace, we can accomplish the following:
\begin{equation}\begin{aligned}
    \kappa(B) &= \int_0^{t^*} tr\left(\gamma^\top e^{t A}B B^\top e^{t A^\top} \gamma\right) dt\\
    & = \int_0^{t^*} tr\left(B^\top e^{t A^\top} \gamma\gamma^\top e^{t A}B\right)dt \\
    &= tr\left(B^\top \int_0^{t^*}e^{t A^\top} \gamma\gamma^\top e^{t A}dt B\right)\\
    &= tr\left(B^\top \Phi B\right).
    \end{aligned}
\end{equation}
Since $B_a^*$ is chosen as eigenvectors of $\Phi$ associated with the same eigenvalue $\lambda$, we have $tr(B^\top \Phi B) = tr(B^\top \lambda B) = \lambda tr(B^\top B) = \lambda m$. Since our goal is to make this denominator as large as possible, it is clear that the global minimum of control energy (global maximum of the denominator) is achieved when $\lambda$ is as large as possible. Therefore, we have $\lambda^* = \sup \Lambda(\Phi) = \lambda_{\max}$.
\end{proof}
This yields an immediate corollary on the energy cost bound:\begin{corollary}\label{cor:flux_energy}
The minimum energy cost is determined by the largest eigenvalue $\lambda^*$ of $\Phi(\gamma)$ as:\begin{equation}
    \bar{\mathcal{E}}_{t^*}^\eta = \frac{\left( \gamma^\top z - \eta\right)^2}{m\lambda^*}.
\end{equation}
\end{corollary} \begin{proof}
The proof follows directly from that of Corollary \ref{cor:mean_energy}. Since $\bar{\mathcal{E}}_{t^*}^\eta = \frac{\alpha(\eta,t^*)^2}{f(B_a^*)}$, and $f(B^*_a) = m\lambda^*$, the solution to optimal energy for control with mean constraints is:\begin{equation}
    \bar{\mathcal{E}}_{t^*}^\eta = \frac{\left( \gamma^\top e^{t^* A}x_0 - \eta\right)^2}{m\lambda^*}.
\end{equation}
\end{proof}

\begin{figure}[h]
    \centering
    \subcaptionbox{Test system}{\includegraphics[width=0.48\linewidth]{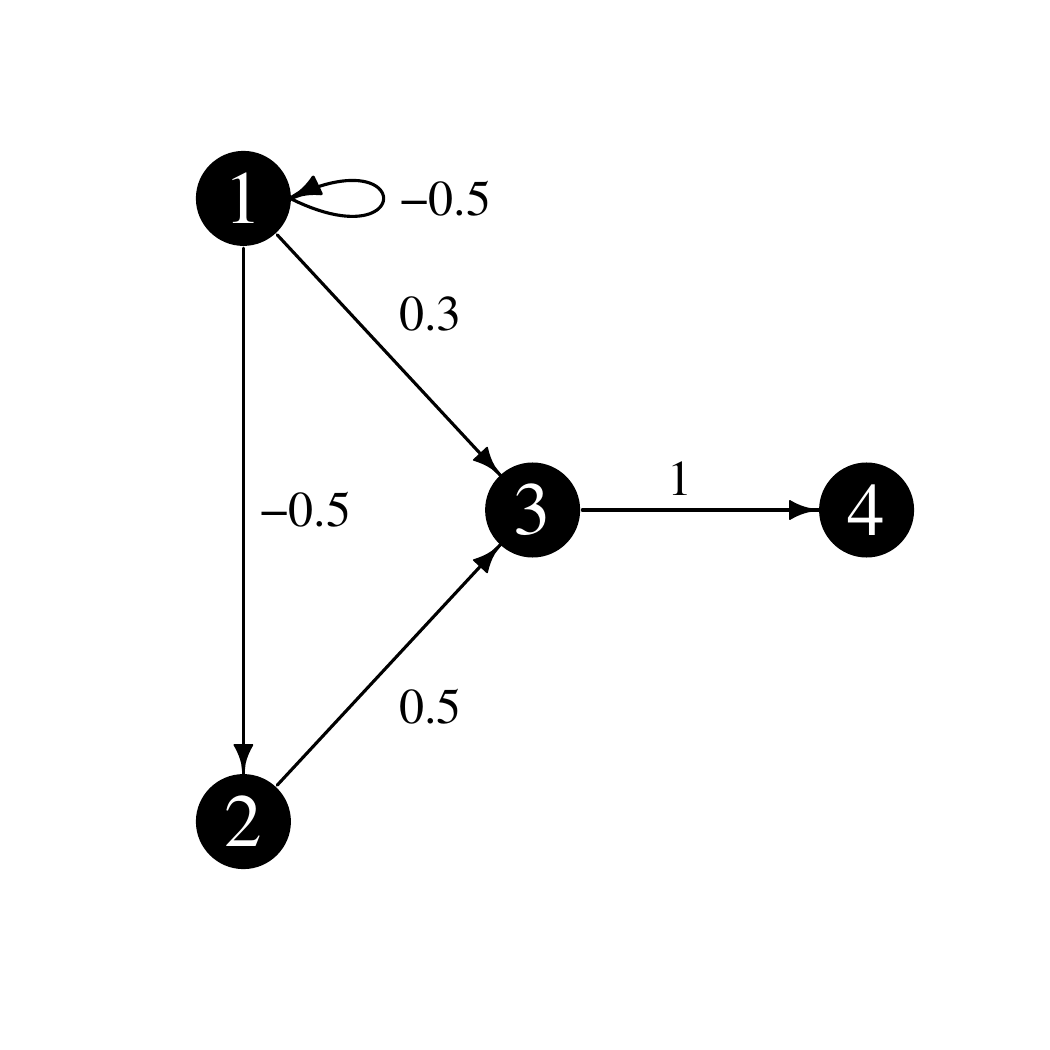}}
    \subcaptionbox{Augmented system}{\includegraphics[width=0.48\linewidth]{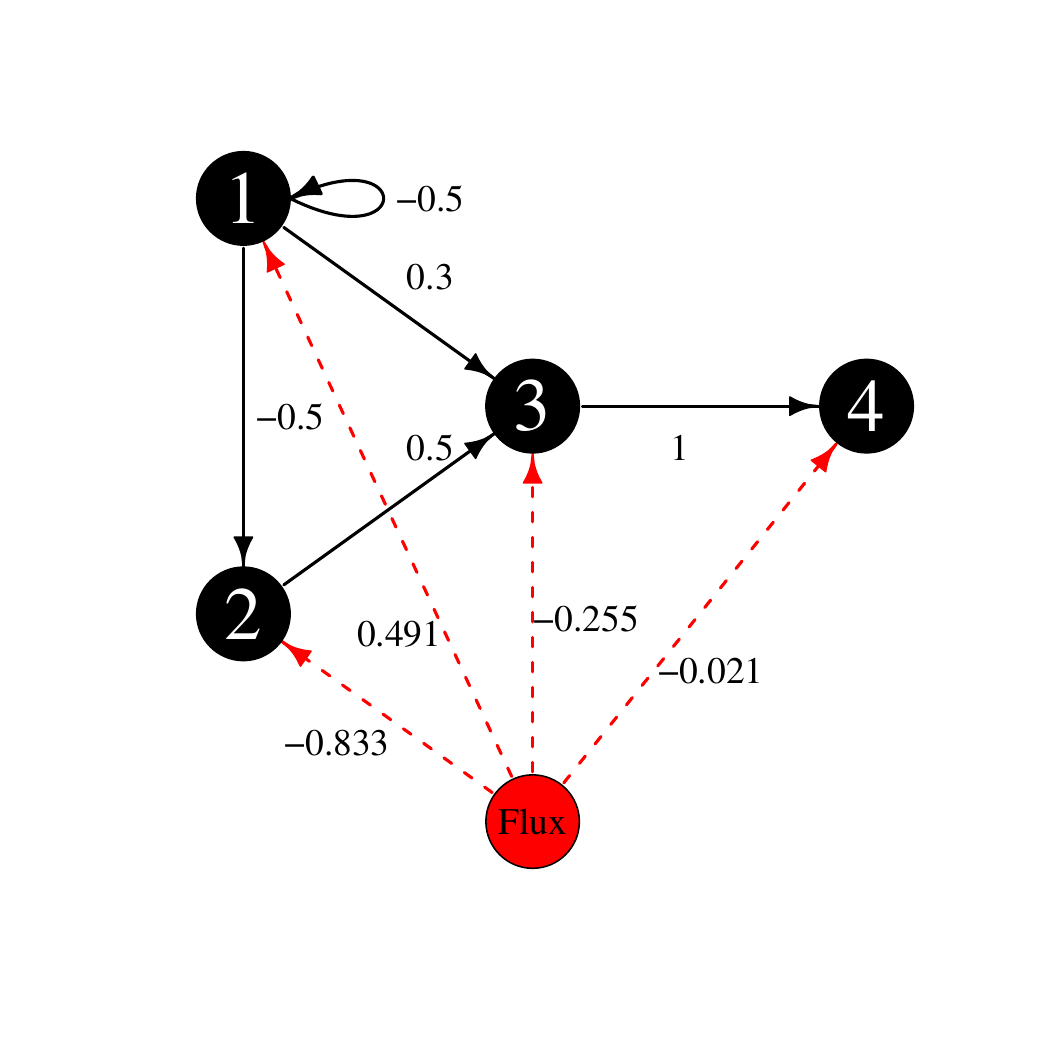}}
    \subcaptionbox{Autonomous dynamics}{\includegraphics[width=0.48\linewidth]{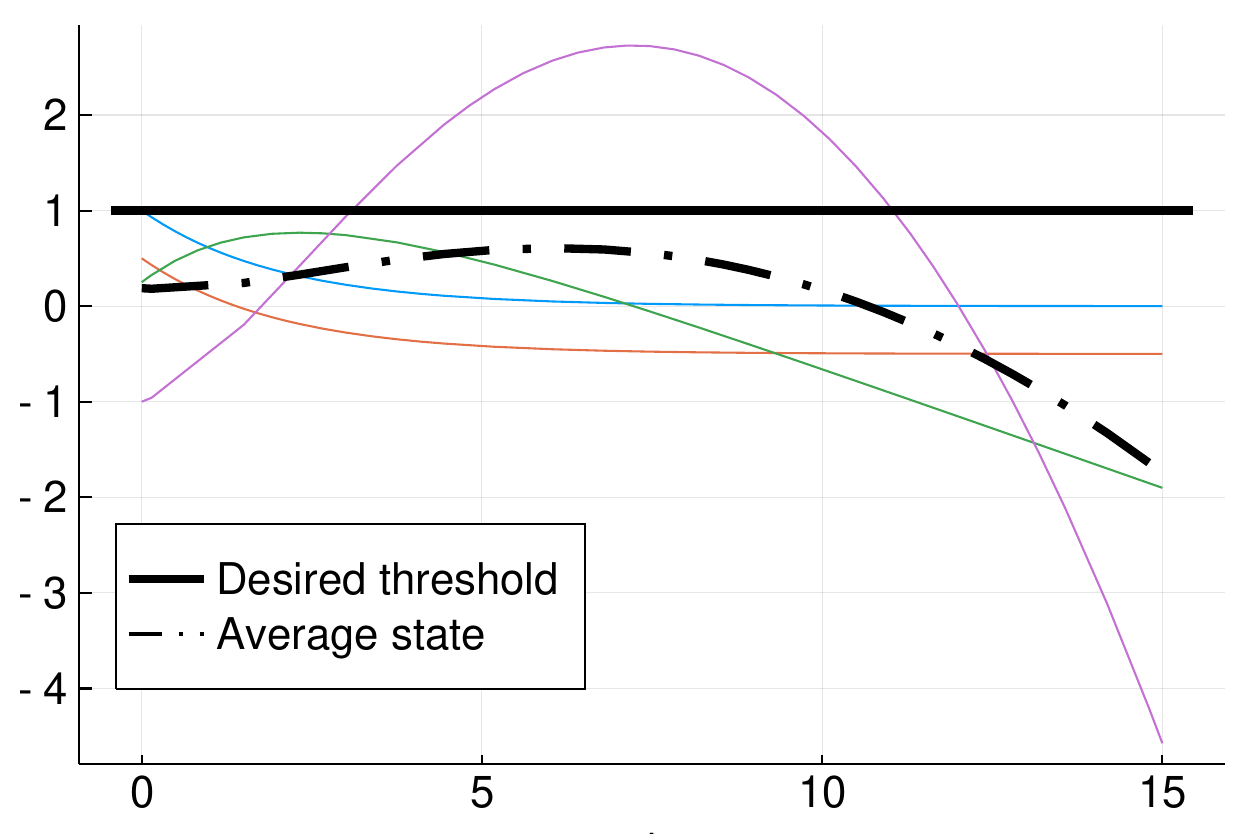}}
    \subcaptionbox{Dynamics under the OMAP}{\includegraphics[width=0.48\linewidth]{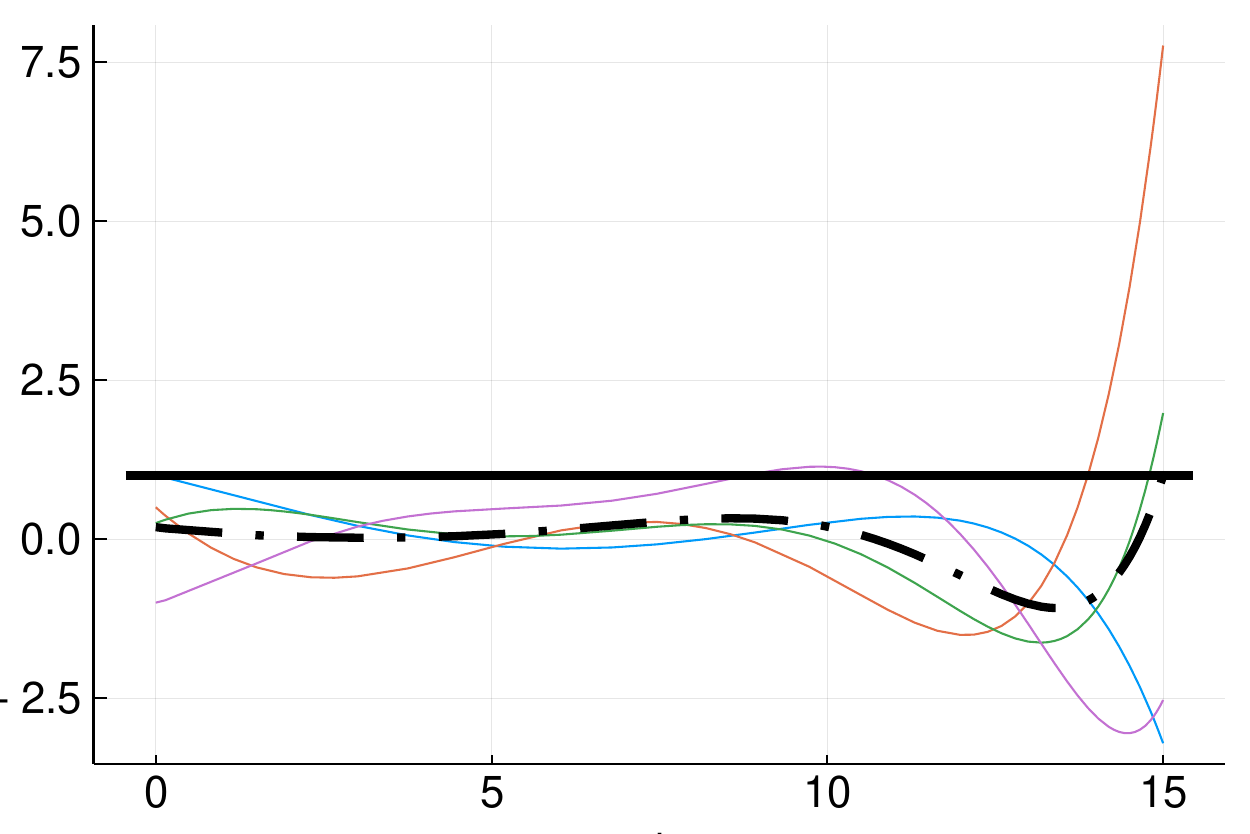}}
    \caption{Dynamics of a test system subject to mean control. Panel (a) gives a graphical depiction of the autonomous system. Panel (b) shows the addition of a control node representing the OMAP and its weights to existing nodes. Panels (b) and (c) show the dynamics of the system state both (c) autonomously and (d) under average control by the OMAP, with the goal $\eta=1$ on $0\leq t \leq 15$.}
\end{figure}

\section{Controlling second moments}
From here, it is natural to consider the problem of controlling higher-order moments of the network state distribution. Returning to the example of platform design for the influence of a social network, second moments may be of particular interest when considering a platform that wishes to strategically influence nodes in order to drive discord or disarray in opinions. On the other hand, such a platform could be designed with an eye toward regulation of discord, in which case the goal of the designer is to optimally place inputs in order to reduce the variance of the state vector faster than the network's autonomous dynamics.
\subsection{State selection on an ellipsoid surface}
Beginning with the problem in its most general form, a constraint on the variance or spread of network states about some point $d\in \mathbb{R}^n$ could be written in the following form:\begin{eqnarray}
    \lVert Ox-d \rVert ^2 = \eta\label{const:variance}.
\end{eqnarray}
This could represent either a prevention of convergence (repulsion or discord) problem, or an aiding of convergence (attraction or regulation), depending on the difference of the threshold $\eta$ with the mean under autonomous evolution.

A particular special case of the variance control problem is that in which the vector $d=\delta \mathds{1}$, for some $\delta\in \mathbb{R}$, such as $\delta = \frac{1}{n}\sum_{i=1}^n z_i$. The constraint (\ref{const:variance}) can be rewritten in vector form as:\begin{equation}
     \left( Ox - d \right)^\top \left( Ox - d \right) = \eta.
 \end{equation}
 Take $L$ as the Cholesky decomposition of the inverse Gramian, where $L$ is a lower triangular matrix, which is full-rank for any controllable system:\begin{equation}W^{-1} = LL^\top. \end{equation}
 Then we can write the state-selection problem as: \begin{equation}\begin{aligned}
    &\underset{x\in \mathbb{R}^{n}}{\text{minimize}} & &\mathcal{E}_{t^*}(x)\\&\text{subject to}&&(Ox-d)^\top (Ox - d) = \eta
    \end{aligned}\end{equation}
    Note that, using Cholesky decomposition from above, defining a constant vector $b\triangleq L^\top z$, this optimization problem can be written as the following:\begin{equation}
        \begin{aligned}
    &\underset{x\in \mathbb{R}^{n}}{\text{minimize}} & &\norm{L^\top x - b}^2 \\&\text{subject to}&&\norm{Ox - d}^2 = \eta
    \end{aligned}\label{opt:variance}
    \end{equation}
    where $O$ is an $n\times n$ matrix. This optimization problem is a case of least-squares with a quadratic constraint \cite{gander1980least,golub1991quadratically}. We can think of this optimization problem as equivalent to finding the stationary points of the energy objective on the surface of the ellipsoid defined by $\norm{Ox-d}^2 = \eta$. The Lagrangian for this optimization problem is:\begin{equation}
        \mathcal{L} = \left(x-z\right)W^{-1}\left(x-z\right) - \lambda \left(\eta - (Ox-d)^\top(Ox-d)\right)\label{lag:var1}.
    \end{equation} The solution to such a problem lies among the solutions to its normal equations, derived from the first-order conditions for $x$ and $\lambda$. In this case, the normal equations are given by:\begin{eqnarray}
        (W^{-1} + \lambda O^\top O )x = W^{-1}z + \lambda O^\top d\label{dldx}\\
        (Ox-d)^\top (Ox-d) = \eta\label{dldl}
    \end{eqnarray}
    If $(W^{-1} + \lambda O^\top O)$ is invertible, then solving the first normal equation (\ref{dldx}) gives:\begin{equation}
        x(\lambda) = \left( W^{-1} + \lambda O^\top O \right)^{-1} \left( W^{-1} z - \lambda d \right).
    \end{equation}
    Finally, plugging this solution for $x$ into the second normal equation (\ref{dldl}) gives the distance function $f(\lambda)$:\begin{equation}
      f(\lambda) = (Ox(\lambda) - d)^\top (Ox(\lambda) - d).
    \end{equation}
    And the \textit{secular equation} of this problem:\begin{equation}
        f(\lambda) = \eta.
    \end{equation}
    The general problem is solved by the smallest nonnegative solution of these normal equations.
    
    \subsection{State repulsion or attraction}
        
    The first special case we will focus on is one in which the controller desires to prevent the dynamics from reaching a given state. In this case, the control signal pushes the state of the system, at the threshold time, further from (or pulls it closer to) a certain given state in the state space. Further, we discuss the special case in which this state is $z$, which would be the final state of the autonomous dynamics; this corresponds to steering the system toward or away from its natural dynamics.
    
    Taking $O=I$ and $d=z$ and have the following result for state-selection:\begin{theorem}
    Given inputs such that $W(B)$ is full-rank, the optimal state which meets the constraint (\ref{const:variance}) at minimum energy cost under the control input schematic $B$ is given by $x^\ast = z + \omega$, where $\omega$ is the leading eigenvector of $W(B)$, normalized to length $\sqrt{\eta}$.
    \end{theorem}\begin{proof}
    First, we will note that $b = L^\top z$ in this case, so the optimization objective from (\ref{opt:variance}) is given by $\norm{L^\top (x-z)}^2$. In the case that $d=z$, we can use a change of variables $x' = z + \omega$ to rewrite (\ref{opt:variance}) as:\begin{equation}
        \begin{aligned}
    &\underset{x'\in \mathbb{R}^{n}}{\text{minimize}} & &\norm{L^\top x' }^2 \\&\text{subject to}&&\norm{x'}^2 = \eta
    \end{aligned}\label{opt:variance2}
    \end{equation}
    The Lagrangian (\ref{lag:var1}) for the repulsion problem then simplifies to:\begin{equation}
        \mathcal{L} = x'^\top W^{-1} x' - \lambda (x'^\top x' - \eta).
    \end{equation}
    First order conditions of this Lagrangian yield the following normal equations:\begin{eqnarray}
    W^{-1} x' = \lambda x' \\
    x'^\top x' = \eta
    \end{eqnarray}
    These normal equations are solved when $\lambda$ is in the spectrum of $W^{-1}$ (and thus $\frac{1}{\lambda} \in \Lambda(W)$), and $x' = \omega$ is its associated eigenvector, with length normalized to $\eta$. Since we have $x' = x - z$, this gives the optimal state $x^\ast$ as $x^\ast = z + \omega$.
    \end{proof}
    The energy cost associated with this state is then straightforward to derive:\begin{corollary}\label{cor:var_energy}
    The minimum energy cost associated with satisfying the constraint (\ref{const:variance}), under control input $B$ is given by \begin{equation}\mathcal{E}(x^\ast ; B) = \sqrt{\eta} \lambda_{\min}(W(B)^{-1}) = \sqrt{\eta} \lambda^{-1}_{\max}(W(B)).\end{equation}
    \end{corollary}\begin{proof}
    Let $\omega = x^\ast - z$ be an arbitrary eigenvector of $W^{-1}$, associated with the eigenvalue $\lambda$. Then the energy cost associated with $x^\ast$ is $\mathcal{E}(x^\ast;B) = \omega^\top W^{-1} \omega = \omega^\top \lambda \omega = \lambda \omega^\top \omega = \lambda \eta$. Since we aim to minimize this energy cost, it is clear that the minimum energy will come when $\lambda = \lambda_{\min} (W(B)^{-1}) = \lambda^{-1}_{\max} (W(B))$.
    \end{proof}
    
    \subsection{Input placement}
With nonlinear outputs, we can still obtain an similar result to Theorem \ref{onecontroller}. In this case, the threshold must be sufficiently large. Formally:\begin{theorem}\label{thm:onecontroller_rep}
    The second-order moment statistic $(x - d)^\top (x- d)$ can be controlled to by a single control input to any threshold $\eta$, provided that \label{mineta}\begin{equation}\eta \geq \norm{d}_2^2 - \norm{d }_\infty^2.\end{equation}
    \end{theorem}
  \begin{proof}
    We wish to prove that for some $\eta$, there exist a state $x\in \text{span}( C(B))$ with $(x-d)^\top (x-d) = \eta$. Let $m=1$, so that $B$ is a vector, and we have $(B-d)^\top (B-d) = B^\top B -2 B^\top d + d^\top d$ For an arbitrary constant $b$, we have $bB \in \text{span}(C(B))$. 
    
    Recall that, since $B$ is a vector, $B^\top B$ is a scalar, so $B^\top B = c \in \mathbb{R}$. Then, we have $bB \in \text{span}(C(B))$, and we seek $b$ such that $(bB-d)^\top (bB-d) = cb^2  - 2b B^\top d + d^\top d = \eta$. Thus, the threshold is reachable if we can find a constant $b$ such that $c b^2 - 2 (B^\top d) b +  (d^\top d - \eta) = 0$. This is a quadratic, with roots at:\begin{equation}
        b^* = \frac{1}{c}B^\top d \pm \frac{1}{c}\sqrt{(B^\top d)^2 - cd^\top d + c\eta}.
    \end{equation}
    Thus, such a $b^* \in \mathbb{R}$ exists if and only if the discriminant is nonnegative, and there exists a vector $B$ such that $(B^\top d)^2 + c\eta \geq c d^\top d$.
    Choose $k$ so that $d_k = \max_i (d_i) = \norm{d}_\infty$, and set $B = \sqrt{c}e_k$, where $e_k$ is a vector of $0$'s with a $1$ in the $k^{th}$ position. This maximizes $B^\top d$, and we have $(B^\top d)^2 \leq c d_k^2$. Thus, a sufficient condition for the existence of a $b^*$ is that $\eta \geq \sum_{i=1}^n d_i^2 - d_k^2$.
    \end{proof}
    While it is very easy to expand the second moment about a given vector, making the threshold smaller makes the problem far more difficult. This corresponds with an inability to ``focus'' network states arbitrarily close to a specific point. This has profound implications -- while it is easy to avoid any given network states, it is more challenging to aid convergence.\footnote{Of course, if the minimum threshold $\eta$ is $0$ for any state $d\in \mathcal{X}$, then we can always force the error from any state to $0$ and thus the system is controllable in the Kalman/PBH sense.} 
    
\subsection{Controlling variance; regulation and discord}

Another important special case of the second moment restriction is to satisfy goals regarding the variance of network states. This is closely related to the ``discord'' problem in social systems studied, for example, by \cite{gaitonde2020adversarial}. In this case, a platform wishes to design inputs that can induce or regulate social discord at minimum cost. The sample variance of states in the state vector $x$ is given by:\begin{equation}
    \bar{\sigma}^2_x = \sum_{i=1}^n \left( x_i - \frac{\mathds{1}^\top x}{n}\right)^2 = \left\lVert x - \frac{\mathds{1}^\top x}{n} \mathds{1} \right\rVert^2.
\end{equation}
    First, we note that we can rearrange this equation to write $\bar{\sigma}^2_x = \lVert Dx \rVert^2$ with a ``centering'' matrix $D \triangleq (I - \frac{1}{n}J)$, (where $J=\mathds{1}\mathds
{1}^\top$ is an $n\times n$ matrix of ones, and $I$ is the identity). Notably, it is clear that $D^k = D$ for all $k>0$; subtracting the mean state will make the mean $0$, so repeated application of the matrix will not change the already centered state. Since, $D$ is also symmetric, we have that $D^\top D = D$. Thus, for the discord state-selection problem\begin{equation}\begin{aligned}
    &\underset{x\in \mathcal{X}^n}{\text{minimize}} & &\mathcal{E}_{t^*}(x;B)\\&\text{subject to}&&\lVert Dx \rVert^2 = \eta
    \end{aligned}\label{meaninputselection3}
\end{equation}
the normal equations are:\begin{eqnarray}
(W^{-1} - \psi D)x = W^{-1}\label{normal1} z\\
x^\top D x = \eta \label{normal2}
\end{eqnarray}
where $\psi$ is the dual variable. To find the solution, we make use of the generalized eigenvalue problem. That is, $\lambda$ is said to be a generalized eigenvalue of $W^{-1}$ and $D$, denoted as $\lambda \in \Lambda(W^{-1},D)$ if the following holds:\begin{equation}
    W^{-1} w = \lambda D w\label{generalizedeigen}
\end{equation}
for some vector $w$. The optimal solution for the Lagrange multiplier is bounded above by the smallest nonnegative generalized eigenvalue of $W^{-1}$ and $D$, which will be denoted as $\lambda_{\min}^+$.

If there is no solution to the normal equations for $\psi$ in the open interval $(0, \lambda_{\min}^+)$, then the solution is $(x^*,\psi^*) = (z-\rho \omega ,\lambda_{\min}^+)$, with $\rho(\eta)$ a real number chosen such that $\lVert z - \rho \omega \rVert^2 = \eta$, and  $\omega$ as the generalized eigenvector associated with $\lambda_{\min}^+$ and length normalized to $\lVert \omega \rVert = 1$.

\begin{figure}[h]
    \centering\hfill
    \subcaptionbox{Network}{\includegraphics[width=0.35\linewidth]{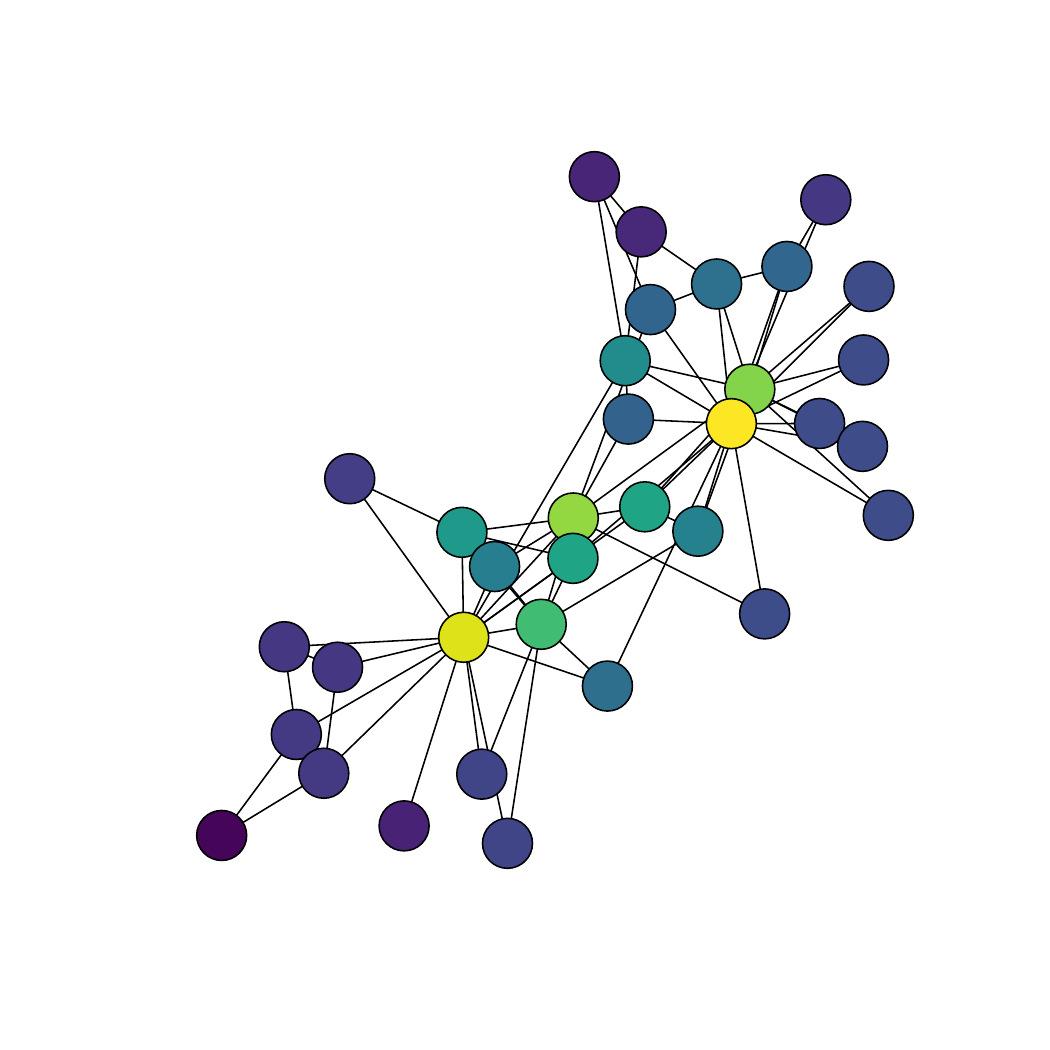}}\hfill
    \subcaptionbox{Autonomous dynamics}{\includegraphics[width=0.48\linewidth]{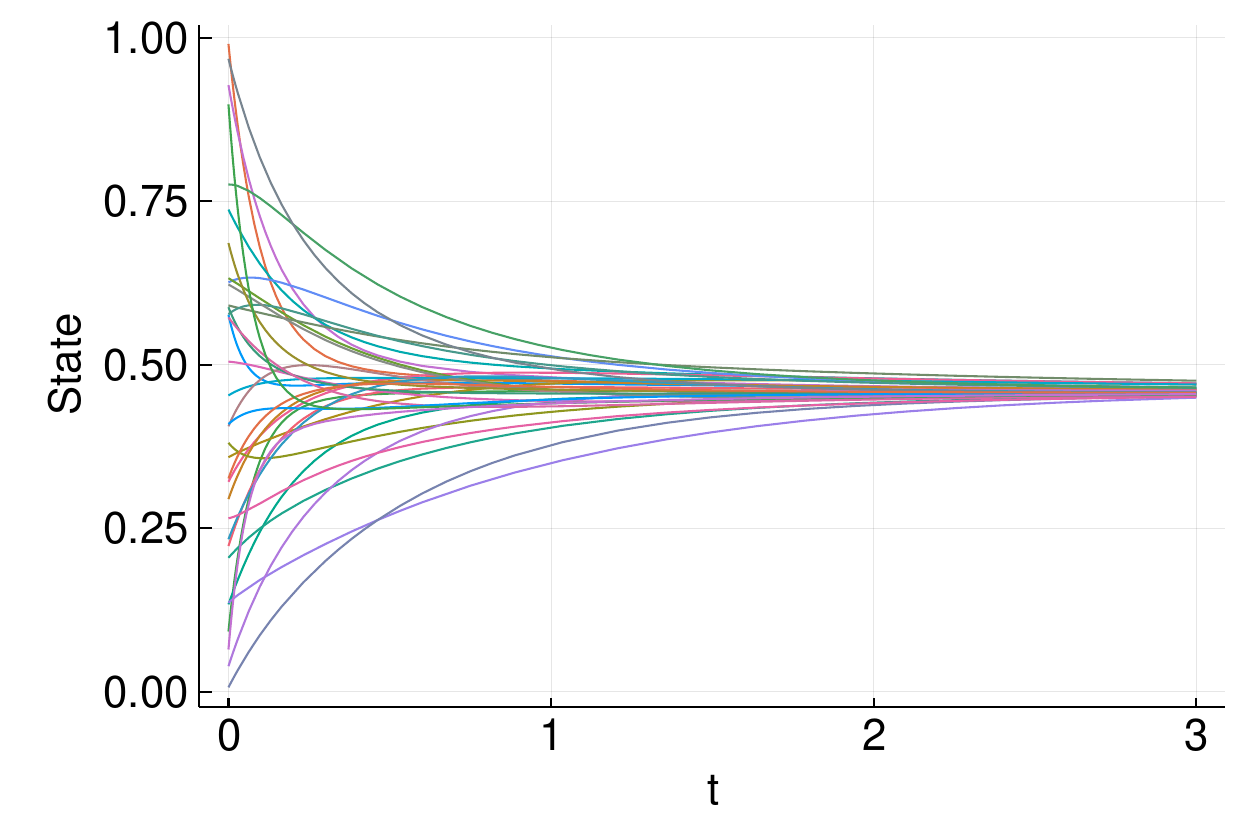}}
    \subcaptionbox{Dynamics under the OMAP}{\includegraphics[width=0.48\linewidth]{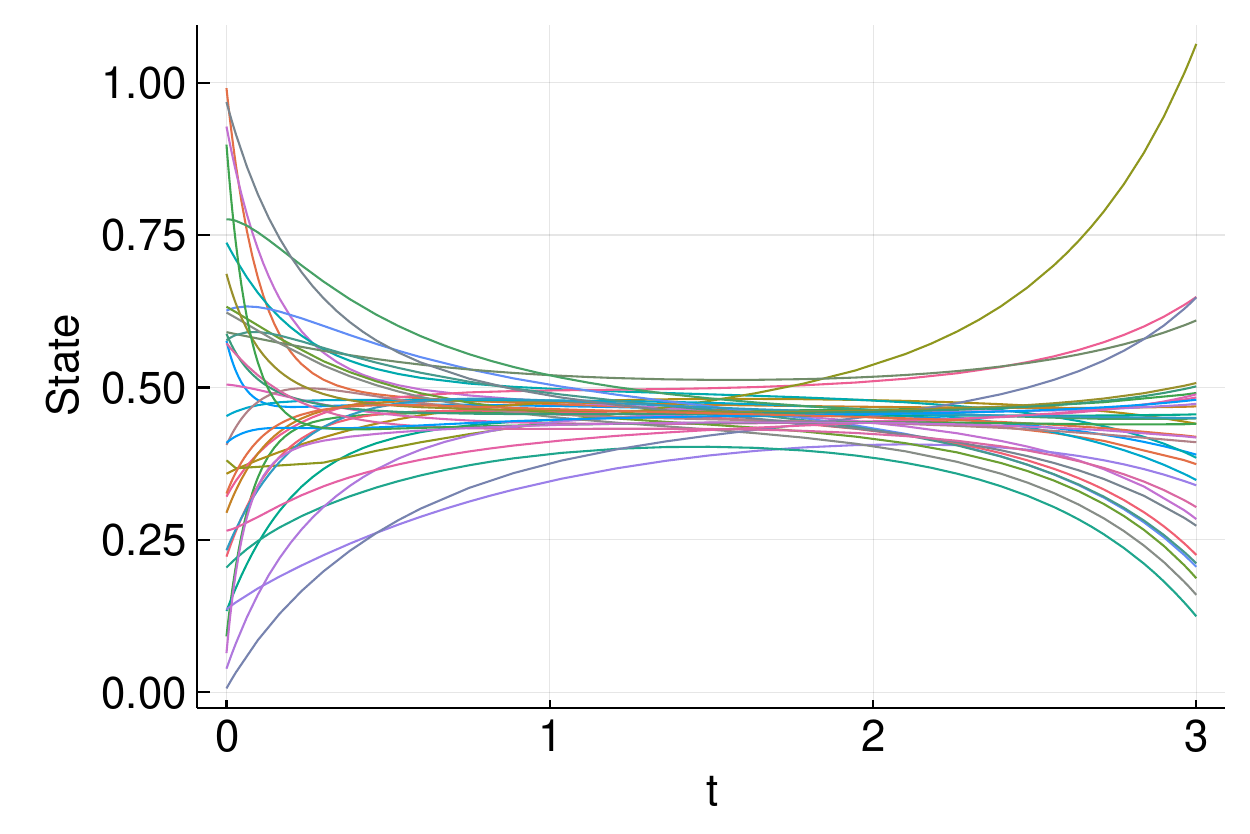}}
    \subcaptionbox{Energy comparison}{\includegraphics[width=0.48\linewidth]{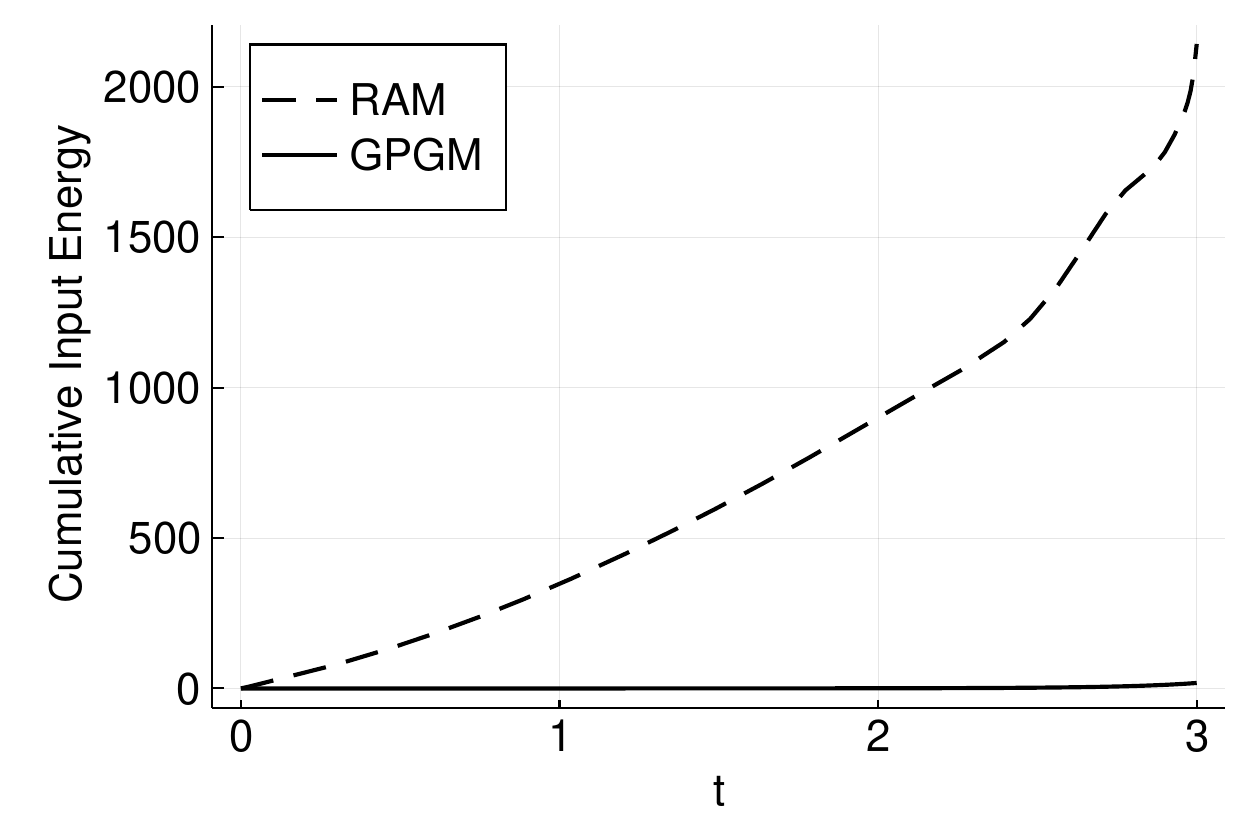}}
    \caption{Dynamics of Zachary's karate club network subject to optimal variance control. Panel (a) shows the Karate club network with nodes colored by eigenvector centrality. Panel (b) shows the Laplacian dynamics of the autonomous system. Panel (c) shows the dynamics under variance control by the OMAP, with ten control nodes and the goal $\eta=1$ on $0\leq t \leq 3$. Panel (d) compares energy costs, of the OMAP with a random input set, and shows that the OMAP produces energy costs that are several orders of magnitude lower.}
\end{figure}

This can be used to obtain an upper bound on the minimum energy cost associated with controlling the variance of the network state vector to a given threshold. Denote the energy cost associated with controlling the variance of network state vector to threshold $\eta$ as $\mathcal{E}_v^*(B,\eta)$.

\begin{theorem}\label{thm:varbound}
The minimum energy cost associated with driving the network to a state vector with variance $\eta$ has the following upper bound:\begin{equation}
    \mathcal{E}_v^* (B,\eta;z) \leq ( \norm{Dz} + \sqrt{\eta} )^2  \lambda_{\min}^+ .\label{eqn:varbound}
\end{equation}
\end{theorem}

\begin{proof}
It is known \cite{gander1980least,gander1989constrained,golub1991quadratically} that there is a solution to the normal equations at the point $(z-\rho(\eta)\omega,\lambda_{\min}^+)$, where $\lambda_{\min}^+$ is the smallest nonnegative generalized eigenvalue of $W^{-1}$ and $D$. Further, this solution is an upper bound on the true solution since there may or may not exist a solution with smaller $\psi$. Denote the positive, scalar-valued function $\rho(\omega,\eta;z)$ as a solution for $\rho$ that satisfies the following secular equation:\begin{equation}
    \norm{D(z-\rho\omega)}^2 = (z-\rho\omega)^\top D (z-\rho\omega) = \eta.
\end{equation} 
We can obtain that $\rho(\omega,\eta;z)^2 \leq \frac{(\norm{Dz}+\sqrt{\eta})^2}{\norm{D\omega}^2}$ by using $\norm{D(z-\rho\omega)}^2=\norm{D(\rho\omega-z)}^2$, and exploiting the reverse triangle inequality:\begin{align*}
    \norm{\rho D\omega} - \norm{Dz}  &\leq \norm{ Dz - \rho D \omega}\\
    \norm{\rho(\omega,\eta;z) D\omega} - \norm{Dz}  &\leq \sqrt{\eta}\\
    \rho(\omega,\eta;z) \norm{D\omega} &\leq \sqrt{\eta} + \norm{Dz}\\
    \rho(\omega,\eta;z) &\leq \frac{\sqrt{\eta}+\norm{Dz}}{\norm{D\omega}}\\
    \rho(\omega,\eta;z)^2 &\leq \frac{(\sqrt{\eta}+\norm{Dz})^2}{\norm{D\omega}^2}.
\end{align*}

It is also known that a higher value of the Lagrange multiplier $\psi$ is associated with a larger objective value. We make use of the Rayleigh quotient $\lambda_{\min}^+ = \frac{\omega^\top M \omega}{\omega^\top D \omega} $ to obtain:\begin{align*}
    \mathcal{E}_v^* (B,\eta;z) &\leq \mathcal{E}(z-\rho(\omega,\eta;z)\omega)\\
    &= (\rho(\omega,\eta;z)\omega)^\top W^{-1} (\rho(\omega,\eta;z)\omega)\\
    &= \rho(\omega,\eta;z)^2 \omega^\top W^{-1} \omega\\
    &= \rho(\omega,\eta;z)^2 \lambda_{\min}^+ \omega^\top D \omega \\
    &\leq \frac{( \norm{Dz} + \sqrt{\eta} )^2}{\norm{D\omega}^2} \lambda_{\min}^+ \norm{D\omega}^2\\
    &= ( \norm{Dz} + \sqrt{\eta} )^2 \lambda_{\min}^+.
\end{align*}
\end{proof}

\section{Input placement with nonlinear outputs}

For the more complex problems of variance constraints, we will proceed by adapting the Projected Gradient Method (PGM) proposed by \cite{li2015minimum}, and replacing the objective with the energy to a minimum cost constrained state rather than a fixed one. To account for the use of various objectives which often need to be iterated to a solution, we suggest a Generalized PGM (GPGM) procedure which nests state selection into the optimization over input configurations.

Following \cite{li2015minimum}, we can project the gradient step direction onto the tangent hyperplane of the constraint sphere to obtain a locally optimal step direction, and then truncate any remaining component of the step direction that would move away from the constraint sphere.\begin{equation}
    N(B) \triangleq (tr (B^\top B) - m)^2
\end{equation}

Differentiating this function with respect to the elements of $B$ gives:\begin{equation}
    \nabla_B N(B) = 2(tr(B^\top B) - m) B
\end{equation}
This gradient is positive as long as $tr(B^\top B) \geq m$, because efficiency is always increasing with the magnitude of an element of $B$. Thus, in order to mitigate an increase in efficiency that comes primarily from increases in $N(B)$, we can project the gradient of the energy functional onto the tangent space of the gradient of this norm function, essentially finding a search direction which impacts energy in a way that does not substantially increase the outgoing connection weights from the input matrix. This is accomplished using the projection operator:\begin{equation}
    \Pi_N (v) = (I - \{\nabla_B N(B) \} \{\nabla_B N(B) \}^\dagger) v
\end{equation}
where $\{ M\}$ represents the vector created by stacking the columns of $M$. This operator projects an arbitrary vector $v$ onto the affine tangent space of the norm function. Since the gradient of the norm function is zero when $tr(B^\top B) = m$, we introduce a small positive constant $\epsilon$ to slightly shrink the radius of the hypersphere to a level at which the gradient of the norm is nonzero.

Because the hypersphere surface $\text{tr}(B^\top B) = m + \epsilon$ is smooth, normalization onto the sphere surface can be done simply, using the following operator:\begin{equation}
 \Pi_s(B) = \left(\sqrt{\frac{m+\epsilon}{\text{tr}(B^\top B)}}\right) B\end{equation}

 The constrained energy function $\widetilde{\mathcal{E}}_c$ is used to represent a function which solves for the minimum energy that can be obtained in the constrained state space $\mathcal{C}$. This could alternatively be obtained by an upper bound such as the one given by Theorem \ref{thm:varbound}. Algorithm \ref{alg1} also depends on a differential operator denoted by $\widetilde{\nabla}$ and a small positive convergence threshold $\delta^*$ with $0<\delta^* \ll 1$.
 
 These algorithms are similar to the PGM and NPGM \cite{li2015minimum,gao2020optimal}, but incorporate different energy functions -- nesting a state optimization to satisfy the metric constraint at the lowest cost. The PGM \cite{li2015minimum} minimizes control energy to a given state subject to the trace constraint, and the NPGM introduced flexibility into this problem by generalizing away from the trace constraint. In a similar fashion, we generalize away from the constraint that the final state is known and exogenous parameter of the problem.
 
 \bigskip \begin{algorithm2e}
\SetKwInOut{Input}{Input}\SetKwInOut{Output}{Output}
\SetAlgoLined
\Input{$\delta^* \in \mathbb{R}$, $\sigma \in \mathbb{R}$, $B' \in \mathbb{R}^{n\times m}$, $A\in \mathbb{R}^{n\times n}$, $x_0\in \mathbb{R}^n$, $\eta \in \mathbb{R}$}
\Output{$B^*_a \in \mathbb{R}^{n\times m}$}\bigskip
$B_k \leftarrow B'$\;
$\delta \leftarrow -\infty$\;
\While{$1-\delta > \delta^*$}{
    $B_{k-1} \leftarrow B_k$\;
    $\{D\} \leftarrow  \{B_{k-1}\} - \sigma \Pi_N (\widetilde{\nabla}_{\{B_{k-1}\}} \widetilde{\mathcal{E}^*}_c (B_{k-1})$)\;
    $\{B_k\} \leftarrow \{\Pi_s \left( D\right)\}$\;
    $\delta \leftarrow \frac{1}{m^2}\{B_k\}^\top \{B_k\}$\;
}
\Return{$B_a^* \leftarrow B_{k}$}\;
\caption{Generalized projected gradient method (GPGM)\label{alg1}}
\end{algorithm2e}\bigskip
 
In cases when the constrained metric is nonlinear, we can use coordinate descent to approximate a local solution to the optimal state; using the iterative procedure to find local improvements in the control energy.

\section{Conclusions}
We have addressed the problem of optimal control for a network system subject to constraints on the distribution of states. When the measure of interest is linear, such as the mean, the system can always be controlled by a single control input. It is thus easy to derive a closed form ranking of nodes that shows how important they are to efficient control of the mean state of the system, based on the spectrum of the observability Gramian. When the measure of interest is nonlinear, the problem must be solved numerically. Notably, it is no longer guaranteed that the measure can be completely controlled with a single input signal. The Generalized Projected Gradient Method (GPGM) uses iterative projections and coordinate-descent to approximate a local solution to the OMAP when the output is nonlinear. In the future, these could be extended to driver node selection strategies through further extensions to the projected gradient method \cite{gao2020optimal}. 

We hope that this work will provide a basis for future study of efficient target or output control of linear dynamical systems with nonlinear outputs. This is a first step in extending PGM-type algorithms to partially nonlinear systems, but there remains a need to extend these results to the fully nonlinear environment.

\section*{Acknowledgements}
We would like to thank James Siderius, Max Gunzburger, Paul Beaumont, and Matt Gentry, along with editors and anonymous referees, for their helpful comments and advice.
\bibliographystyle{IEEEtran}
\bibliography{main_v4}

\end{document}